\documentclass[a4paper,12pt]{article}
\usepackage{amsmath}
\usepackage{amssymb}
\usepackage{tabularx}
\usepackage{enumerate}

\newtheorem{thm}{Theorem}[section]
\newtheorem{lem}[thm]{Lemma}

\newtheorem{rmk}{Remark}[section]
\newtheorem{pppp}{Proof}

\newenvironment{pf}{\begin{pppp} \em}{\mbox{}\hfill\qed\end{pppp}}
\newcommand{\qed}{\hspace{1em}\mbox{\raisebox{0.65ex}{\fbox{}}}}

\numberwithin{equation}{section}
\newcommand{\ol}{\overline}

\newcommand{\be}{\begin{equation}}
\newcommand{\ee}{\end{equation}}
\newcommand\bes{\begin{eqnarray}} \newcommand\ees{\end{eqnarray}}
\newcommand{\bess}{\begin{eqnarray*}}
\newcommand{\eess}{\end{eqnarray*}}

\begin{document}

\thispagestyle{empty}

\title{{\bf The spreading fronts in a mutualistic
model with delay\thanks{The work is supported by Supported in part by a NSFC Grant No. 11171158,
NSF of Jiangsu Education Committee No. 11KJA110001, "333" Project of Jiangsu Province  Grant No. BRA2011173.}}}
\date{\empty}
\author{Mei Li$^{a,}$ $^b$\footnote{The corresponding author. E-mail address: limei6606@tom.com}\\
{\small $^a$ School of Applied Mathematics, Nanjing University of Finance and Economics,}\\
{\small Nanjing 210023,PR China}\\
{\small $^b$ Institute of mathematics, Nanjing Normal University,}\\
{\small Nanjing 210023,PR China}}
\maketitle

\begin{quote}
\noindent {\bf Abstract.} {\small This article is concerned with a
system of semilinear parabolic equations with two free boundaries
describing the spreading fronts of the invasive species in a mutualistic ecological model. The local existence and
uniqueness of a classical solution are obtained and the asymptotic
behavior of the free boundary problem is studied. Our results indicate that two free boundaries tend monotonically to finite or infinite
at the same time, and the free boundary problem admits a global slow solution with unbounded free boundaries if the
inter-specific competitions are strong, while if the inter-specific
competitions are weak, there exist the blowup solution and global
fast solution.}



\end{quote}

\section{Introduction}

In this paper, we consider the following parabolic
system with moving boundaries:
\begin{eqnarray}
\left\{
\begin{array}{lll}
u_{t}-d_1 u_{xx}=u(a_1-b_{1}u+c_{1}v(t-\tau_1,x)),\; & t>0, \ g(t)<x<h(t),  \\
v_{t}-d_2 v_{xx}=v(a_2+b_{2}u(t-\tau_2,x)-c_{2}v),\; & t>0, \ -\infty<x<\infty, \\
u(t, x)=0,\; & t\geq 0, \ x<g(t)\, \textrm{or}\, h(t)<x,  \\
u=0,\quad h'(t)=-\mu \frac{\partial u}{\partial
x},\quad & t>0, \  x=h(t),  \\
u=0,\quad g'(t)=-\mu \frac{\partial u}{\partial
x},\quad & t>0, \  x=g(t),  \\
-g(0)=h(0)=b, \quad (0<b<\infty),&\\
u(t, x)=u_{0}(x)\geq 0,\; &-b\leq x\leq b,-\tau_{2}\leq t\leq 0,\\
v(t, x)=v_{0}(x)\geq 0,\; &-\infty\leq x\leq \infty,-\tau_1\leq t\leq 0,\\
\end{array} \right.
\label{f1}
\end{eqnarray}
where $x=h(t)$ and $x=g(t)$ are the moving boundaries to be determined. Here
 $a_i$, $b_i$ and $c_{i}$ ($i=1,2$) are positive constants.
$u_0(x)$ and $v_0(x)$ are initial functions satisfying
\begin{eqnarray}
\left\{
\begin{array}{ll}
u_0\in C^1[-b, b]\cap C^2(-b,b),\,  v_0\in L^\infty (-\infty, \infty)\cap C^2(-\infty, \infty),  \\
u(-b)=u(b)=0,\, u'_0(-b)>0,\ u'_0(b)<0,\\
u_0>0,v_0>0\ \textrm{in}\ (-b, b).
\end{array} \right.
\label{a11}
\end{eqnarray}
This system describes the cooperating two-species Lotka-Volterra
model, where the native species $(v)$ migrates in the habitat $(-\infty, \infty)$ and
the invasive species $(u)$ is initially limited in a special part and
disperses through random diffusion only in $g(t)<x<h(t)$. In biological terms, the unknowns $u(x,t)$
and $v(x,t)$ represent the spatial densities of the species at time $t$ and location $x$,
$a_i$ is its respective net birth rate and the constant $d_i>0$ is the
diffusion coefficient. The coefficients $b_1$ and $c_2$ measure the
intra-specific competitions whereas $b_2$ and $c_1$ represent
inter-specific cooperation.

The corresponding problem on a fixed domain  transforms into a Lotka-Volterra mutualistic model:
 \bes \left\{
\begin{array}{ll}
u_t=d_1\Delta u+u(a_1-b_{1}u+c_{1}v)\; &\textrm{for} \;   t>0, \ x\in\Omega, \\
v_t=d_2\Delta v+v(a_2+b_{2}u-c_{2}v)\; &\textrm{for} \;   t>0, \ x\in\Omega,
\end{array} \right.
\label{a2} \ees
which can be interpreted in biological terms that
the presence of one species encourages the growth of the other
species. Pao \cite{P} displied that the solution of (\ref{a2}) under Dirichlet boundary condition with any
initial data is unique and global when $b_2c_1<b_1c_2$, while the
blowup solutions are possible when the two species are strongly
mutualistic ($b_2c_1>b_1c_2$), which means that the geometric mean
of the interaction coefficients exceeds that of population
regulation coefficients.

 The conditions on the free boundaries are
$h'(t)=-\mu u_{x}(t, h(t))$ and $g'(t)=-\mu u_{x} (t,g(t))$,
which are called the Stefan conditions. Here it means that the amount of the species flowing across
the free boundary is increasing with respect to the moving length, see
\cite{Li} in detail, and $\mu$ is positive constant.

  Recently Kim and Lin \cite{KLL} studied  the corresponding system of semilinear parabolic with a free boundary
\begin{eqnarray}
\left\{
\begin{array}{ll}
u_{t}-d_1 u_{xx}=u(a_1-b_{1}u+c_{1}v),\; & t>0, \ 0<x<h(t),  \\
v_{t}-d_2 v_{xx}=v(a_2+b_{2}u-c_{2}v),\; & t>0, \ 0<x<\infty, \\
u(t, x)=0,\; & t\geq 0, \  h(t)<x<\infty,  \\
u=0,\quad h'(t)=-\mu \frac{\partial u}{\partial
x},\quad & t>0, \  x=h(t),  \\
\frac{\partial u} {\partial x}(t,0)=\frac{\partial v} {\partial x}(t,0)=0,\;& t>0,\\
h(0)=b, \quad (0<b<\infty),&\\
u(0, x)=u_{0}(x)\geq 0,\; &0\leq x\leq b,\\
v(0, x)=v_{0}(x)\geq 0,\; &0\leq x\leq \infty,\\
\end{array} \right.
\label{f2}
\end{eqnarray}
the blowup solution and global fast solution are given.

 In the absence of $v$ and the nonlinear reaction term for $u$, problem (\ref{f1})
 is reduced to one phase Stefan problem, which accounts for phase transitions between solid and fluid states
such as the melting of ice in contact with water \cite{Ru}. Stefan
problem has been studied by many authors, see \cite{
Cu, DL, DL2, FR, KN, LZZ, Ol, Ra, Ru,T}.

As to the one-phase Stefan problem for the heat equation with a superlinear reaction term
\begin{eqnarray}
\left\{
\begin{array}{lll}
u_{t}- u_{xx}=u^{1+p},\; & t>0, \ 0<x<h(t), \\
h'(t)=-\frac{\partial u}{\partial
x},\quad &  t>0, \ x=h(t), \\
\frac{\partial u}{\partial x}(t, 0)=u(0, h(t))=0,\; &t>0,\\
u(0, x)=u_{0}(x)\geq 0,\; &0\leq x\leq b,\, \, h(0)=b,
\end{array} \right.
\label{f4}
\end{eqnarray}
it was shown in \cite{SP2, SP1} that  all global solutions are
bounded and decay uniformly to $0$ as $t\rightarrow\infty$ if the
initial data is small, while if the initial date is big, the
solution will blow up in a finite time. Moreover they showed that
there exist global solutions with slow decay and unbounded free
boundary.

The free boundary problems associated with the ecological models have attracted
considerable research attention in the past due to their relevance in applications, see
for example, \cite{HM, Li, MY1, MY2, MY3} and the references therein.

Motivated by Kim and Lin \cite{KLL}, we are interested asymptotic behaviors of the
solution for two free boundaries problem (\ref{f1}), especially a more detailed category
about the global solutions. We will show that if $b_1c_2>b_2c_1$, there exists a global slow solution of (\ref{f1}), while if $b_1c_2<b_2c_1$ there exist a blowup solution and global fast solution of (\ref{f1}).

Throughout this paper, a solution $(u, v, g, h)$ of (\ref{f1}) is said
to be classical if $u\in C([0, T)\times [g(t), h(t)])$ $\bigcap
C^{1,2}((0, T)\times (g(t), h(t))$, $v\in C([0, T)\times (-\infty, \infty ))$
$\bigcap C^{1,2}((0, T)\times (-\infty, \infty))$ $\bigcap C([0, T)\times
L^\infty(-\infty, \infty))$ and $h,g\in C^1[0, T)$ with $T_{\max}\leq +\infty$ and
satisfies (\ref{f1}), where $T_{\max}$ denotes the maximal existing time of solution. If $T_{\max}=+\infty$, we say the solution exists
globally whereas if the solution ceases to exist for some finite
time, that is, $T_{\max}<+\infty$ and $\lim_{t\to T_{\max}} (||u(t,
x)||_{L^\infty([g(t), h(t)])}+||v(t, x)||_{L^\infty (-\infty, \infty)})\to
+\infty $, we say that the solution blows up.  If $T_{\max}=\infty$ and
$h_\infty:=\lim _{t\to \infty}h(t)<\infty, g_\infty:=\lim _{t\to \infty}g(t)>-\infty$,
the solution is called global fast solution since that the solution decays
uniformly to $0$ at an exponential rate, while If $T_{\max}=\infty$ and
$h_\infty=\infty, g_\infty=-\infty$, it is called
global slow solution, whose decay rate is at most polynomial, see \cite{SP2,SP1}.

 We now briefly give an outline of the paper. In Section 2, local existence and
uniqueness of two free boundaries problem (\ref{f1}) are obtained by
using Schauder fixed point theorem. Results pertaining to global solution
 for the case $b_1c_2>b_2c_1$ are presented in
Section 3, and in Section 4, results regarding nonglobal solutions and
global fast solution  for the case $b_1c_2<b_2c_1$ are established.

\section{Local existence and uniqueness}

 We first prove the following local existence and uniqueness
results of the solution to $(\ref{f1})$ by virtue of the Schauder
fixed point theorem:

\begin{thm} There exists a $T>0$ such that problem $(\ref{f1})$
admits a unique solution $$(u, v, h, g)\in
C^{1+\alpha,(1+\alpha)/2}(\overline{D}_{1,T})\times C^{1+\alpha,
(1+\alpha)/2}(\overline{D}_{2,T}) \times C^{1+\alpha/2}[0,T]\times C^{1+\alpha/2}[0,T],$$
furthermore
$$\|u\|_{C^{1+\alpha,(1+\alpha)/2}(\overline{D}_{1,T})}+\|v\|_{C^{1+\alpha,(1+\alpha)/2}(\overline{D}_{1,T})}$$
\begin{eqnarray}
+\|h\|_{C^{1+\alpha/2}[0,T]}+\|g\|_{C^{1+\alpha/2}[0,T]}\leq
C,\label{b12}
\end{eqnarray}
where  $D_{1,T}=(0, T]\times (g(t), h(t))$, $D_{2,T}=(0, T]\times (-\infty,
+\infty)$,  $0<\alpha <1$, and $C,T$ only depend on $b$,
$\|u_0\|_{C^{2}[-b, b]}$ and $\|v_0\|_{C^{2}(-\infty,
+\infty)}$.
\end{thm}
\begin{pf}
As in  \cite{CF1} and \cite{DL}, we first straighten the double free
boundaries. Let $\zeta (y)$ be a function in $C^3(-\infty,+\infty)$
satisfying
$$\zeta (y)=1\quad  \textrm{if} \,\, |y-b|<\frac{b}8,$$
$$\zeta (y)=0\quad  \textrm{if} \,\, |y-b|>\frac b 2,\, |\zeta '(y)|<\frac 6b,$$
$$\xi(y)=\zeta(-y).$$
Let a transformation be
$$(t, x)\rightarrow (t, y), \textrm{where}\,\, x=y+\xi (y)(g(t)+b)+\zeta (y)(h(t)-b),
\quad -\infty< y<\infty.$$ As long as $$\max \{|g(t)+b|,|h(t)-b|\}\leq \frac b8,$$ the
above transformation is a diffeomorphism from $(-\infty, +\infty)$ onto
$(-\infty, +\infty)$. Moreover,  the free boundary $x=h(t), x=g(t)$ becomes
the lines $y=b, y=-b$ respectively. Now, a straightforward computation yields
\begin{eqnarray*}
&\frac {\partial y}{\partial x}=\frac 1{1+\xi'(y)(g(t)+b)+\zeta'(y)(h(t)-b)}\equiv
\sqrt{A(g(t),h(t),y)}\equiv{C(g(t),h(t),y)},&\\
&\frac {\partial^2 y}{\partial x^2}=-\frac {\xi''(y)(g(t)+b)+\zeta
''(y)(h(t)-b)}{[1+\xi'(y)(g(t)+b)+\zeta'(y)(h(t)-b)]^3}\equiv
B(g(t),h(t),y),&\\
&-\frac {\partial y}{\partial t}=\frac
{\xi(y)g'(t)+\zeta(y)h'(t)}{1+\xi'(y)(g(t)+b)+\zeta'(y)(h(t)-b)}\equiv
{C(g(t),h(t),y)}[{\xi(y)g'(t)+\zeta(y)h'(t)}].&
\end{eqnarray*}

If we set
$$u(t, x)=u(t, y+\xi(y)(g(t)+b)+\zeta (y)(h(t)-b))=w(t, y),$$
$$v(t, x)=v(t, y+\xi(y)(g(t)+b)+\zeta (y)(h(t)-b))=z(t, y),$$
then
$$ u_t=w_t-{\xi(y)g'(t)+\zeta(y)h'(t)}C(g(t),h(t),y)w_y,$$
$$ v_t=z_t-{\xi(y)g'(t)+\zeta(y)h'(t)}C(g(t),h(t),y)z_y,$$
$$ u_x={C(g(t),h(t),y)}w_y, \quad  v_x={C(g(t),h(t),y)} z_y,$$
$$ u_{xx}=A(g(t),h(t),y)w_{yy}+B(g(t),h(t),y)w_y, $$
$$ v_{xx}=A(g(t),h(t),y) z_{yy}+B(g(t),h(t),y)z_y$$
and problem (\ref{f1}) turns into
\begin{eqnarray}
\left\{
\begin{array}{lll}
w_{t}-Ad_1 w_{yy}-[Bd_1+({\xi(y)g'(t)+\zeta(y)h'(t)})C]w_y \\
=w(a_1-b_{1}w+c_{1}z(t-\tau_1,y)),\; &t>0, \
-b<y<b, \\
z_{t}-Ad_2 z_{yy}-[Bd_2+({\xi(y)g'(t)+\zeta(y)h'(t)})C]z_y \\
=z(a_2-b_2w(t-\tau_2,y)-c_2z),\; & t>0, \ -\infty<y<\infty, \\
w(t,y)=0,\; &t\geq 0, \ -\infty<y<-b,\\
w(t,y)=0,\; &t\geq 0, \ b<y<\infty,\\
w=0,\quad h'(t)=-\mu\frac{\partial w}{\partial
y},\quad &t>0, \ y=b,\\
w=0,\quad g'(t)=-\mu\frac{\partial w}{\partial
y},\quad &t>0, \ y=-b,\\
h(0)=-g(0)=b, \quad (0<b<\infty).&\\
w(t, y)=u_{0}(y)\geq 0,\; &-b\leq y\leq b,-\tau_2\leq t\leq 0\\
z(t, y)=v_{0}(y)\geq 0,\; &-\infty\leq y<\infty,-\tau_1\leq t\leq 0,
\end{array} \right.
\label{b}
\end{eqnarray}
where  $A=A(g(t), h(t), y)$, $B=B(g(t), h(t), y)$, $C=C(g(t), h(t), y)$, $u_0(y)\in
C^1[-b, b]\cap C^2(-b,b)$ and $v_0(y)\in L^\infty (-\infty, \infty)
\cap C^2(-\infty, \infty)$.

Let $g_1=-\mu u'_0(-b),h_1=-\mu u'_0(b)$, $D^*_{1,T}=(0, T]\times (-b, b)$ and
$0<T<\min(\frac {b}{8(1+g_1)},\frac {b}{8(1+h_1)})$, choosing
 \begin{eqnarray*}
&&\mathcal{D}_1=\{w(t, y)\in
C(\overline{D}^*_{1,T}): \ w(t, y)=u_0(y)\},\\
&&\mathcal{D}_{1T}=\{w\in \mathcal{D}_1: \
\sup_{0\leq t\leq T,\ -b\leq y\leq b}|w(t, y)-u_0(y)|\leq 1\},\\
&&\mathcal{D}_2=\{z(t, y)\in C(\overline{D}_{2,T}): \ z(t, y)=v_0\},\\
&&\mathcal{D}_{2T}=\{z\in \mathcal{D}_2:\  \sup_{0\leq t\leq T,\ -\infty< y<\infty} |z(t, y)-v_0(y)|\leq 1\},\\
&&\mathcal{D}_3=\{g(t)\in C^1[0,T]: \ g(0)=-b,\  g'(0)=g_1\},\\
&&\mathcal{D}_{3T}=\{g(t)\in\mathcal{D}_3: \ \sup_{0\leq t\leq
T}|g'(t)-g_1|\}\leq 1\},\\
&&\mathcal{D}_4=\{h(t)\in C^1[0,T]: \ h(0)=b,\  h'(0)=h_1\},\\
&&\mathcal{D}_{4T}=\{h(t)\in\mathcal{D}_4: \ \sup_{0\leq t\leq
T}|h'(t)-h_1|\}\leq 1\}.
\end{eqnarray*}
 It's well known that $\mathcal
{D}_{1T}\times\mathcal {D}_{2T}\times\mathcal {D}_{3T}\times\mathcal {D}_{4T}$ is a closed
convex set in $C(\overline{D}^*_{1,T})\times
C(\overline{D}_{1,T})\times C^1[0,T]\times C^1[0,T]$.

Next, we can obtain the existence and uniqueness by using the contraction mapping theorem
 as in \cite{CF1, DL} with some obvious adaptation. For brief, we omit it here.
\end{pf}

\begin{thm} The double free boundaries in problem $(\ref{f1})$ are sternly monotone,
namely, for any solution on $[0, T]$,we have
 $$h'(t)>0\ \textrm{and} \ g'(t)<0\ \textrm{for} \ 0\leq t\leq T.$$
\end{thm}
\begin{pf}
Using the Hopf Lemma to the system of$(\ref{f1})$, we  deduce that
 $$u_x(t, h(t))<0,\ u_x(t,g(t))>0 \ \textrm{for} \ 0\leq t\leq T.$$
Then, combining the above two inequalities with the Stefan conditions in $(\ref{f1})$,
the result can be deduced.
\end{pf}

Furthermore, the double free boundaries $g(t)$ and  $h(t)$
have another notable properties which will be showed below.

\begin{thm} Let $(u, v, g, h)$ be a solution of system $(\ref{f1})$ in $[0, T_{\max})\times [g(t), h(t)]$.
Then $g(t)$ and $h(t)$ satisfy
$$-2b<g(t)+h(t)< 2b, \  t\in [0,T_{\max}).$$
\end{thm}

\begin{pf}
It follows from continuity that $g(t)+h(t)<2b$  for small $t>0$. Define
$$\emph{T}:= \sup\{s:\, g(t)+h(t)<2b,\ t\in[0,s)\}.$$
We can deduce that $\emph{T}=\emph{T}_{\max}$ in the following proof by contradiction.
Suppose that $\emph{T}<\emph{T}_{\max}$, we then have
$$g(t)+h(t)<2b,\, t\in[0,\emph{T}),\quad g(\emph{T})+h(\emph{T})=2b.$$
Hence
\begin{equation}
 g'(\emph{T})+h'(\emph{T})\geq 0.
\label{b12}
\end{equation}
In order to obtain a contradiction, we define the function
$\mathcal{F}(t,x):=u(t,x)-u(t,-x+2b)$
on the region
$$\Omega' =\{(t, x): \ 0\leq t\leq T,\  b\leq x\leq h(t)\}.$$
Directly calculating $F$ shows that it satisfies
$$F_{t}=F_{xx}+c(t,x)F,\ 0<t\leq\emph{T},\, b<x<h(t),$$
with some $c(t,x)\in L^\infty (\Omega')$ and
$$F(t, b)=0,\ F(t,h(t))<0,\ 0<t<\emph{T}.$$
Moreover,
$$F(\emph{T},h(\emph{T}))=u(\emph{T},h(\emph{T}))-u(\emph{T},-h(\emph{T})+2b)=u(\emph{T},h(\emph{T}))-u(\emph{T},g(\emph{T}))=0.$$
Then we have
$$F(t,x)<0,\, (t,x)\in (0,\emph{T}]\times (b, h(t)),$$
and
$$F_{x}(\emph{T},h(\emph{T}))<0$$
by applying the strong maximum principle and the Hopf lemma.
However
$$F_x(\emph{T},h(\emph{T}))=u_x(\emph{T},h(\emph{T}))+u_x(\emph{T},g(\emph{T}))=-[g'(\emph{T})+h'(\emph{T})]/\mu,$$
namely
$$g'(\emph{T})+h'(\emph{T})>0,$$
which contradicts to (\ref{b12}). Therefore we claim that
$g(t)+h(t)<2b$ for all $0<t<\emph{T}_{\max}.$
Similarly we can prove $g(t)+h(t)>-2b$ for all $0<t<\emph{T}_{\max}$.
\end{pf}

Theorem 2.1 implies that there exists a $T$ such that the solution exists
in time interval $[0,T]$, and the solution can be further extended to $[0,T_{\max})$  with $T_{\max}\leq +\infty$ by Zorn's lemma.
The maximal exist time of the solution
$T_{\max}$ depends on a prior estimate with respect to
$||u||_{L^\infty}$, $||v||_{L^\infty}$ and $g'(t), h'(t)$. Next we will
give that if $||u||_{L^\infty}<\infty$, the solution is global. For this purpose we first provide the following lemma:
\begin{lem}
Suppose that $\overline M\triangleq||u||_{L^\infty([0, T]\times[g(t), h(t)])}<\infty$. Then the solution of
the free boundary problem $(\ref{f1})$ satisfies
$$0\leq v\leq M_2(\overline M)\ \  \textrm{for} \  0\leq t\leq T,\  -\infty\leq x<\infty,$$
$$ 0<-g'(t)\leq M_3(\overline M) \ \  \textrm{for} \ \  0\leq t\leq T,$$
$$ 0<h'(t)\leq M_4(\overline M) \ \  \textrm{for} \ \  0\leq t\leq T,$$
\end{lem}
where $M_2, M_3$ and $M_4$ are independent of $T$.
\begin{pf}
Since that $v_t-d_2v_{xx}\leq v(a_2+b_{2}\overline M-c_{2}v)$ for
$0<t\leq T$, $-\infty<x<\infty$, the estimate for $v$ is directly from the
Phragman-Lindelof principle.

Set $$\Omega =\{(t, x): \ 0<t\leq T,\  g(t)<x<g(t)+\frac 1M\}$$ and
constitute an auxiliary function
$$w(t, x)=\overline M[2M(x-g(t)-M^2(x-g(t))^2].$$
In the following proof, we will choose $M$ such that $w(t, x)$ is the supersolution of
$u(t, x)$ in $\Omega$.

Tedious but fairly straightforward computation show that
$$w_t=2\overline M M(-g'(t))\big(1-M(x-g(t))\big)\geq 0,$$
$$-w_{xx}=2\overline M M^2,$$
$$u(a_1-b_{1}u+c_{1}v)\leq \overline M(a_1+c_{1}M_2).$$
It follows that $$w_t-d_1w_{xx}\geq \overline M(a_1+c_{1}M_2)\geq
u(a_1-b_{1}u+c_{1}v)$$ if $M^2\geq \frac {a_1+c_{1}M_2}{2d_1}$.
On the other hand,
$$w(t, g(t)+\frac 1M)=\overline M\geq u(t, g(t)+\frac 1M),$$
$$w(t, g(t))=0=u(t, g(t)).$$ Recalling that $u_0(-b)=0$ and
$u'_0(-b)=-g_1/\mu$ gives that there exists $0<\delta <b$ such that
$u_0(x)\leq \frac 34 \overline M$ and $|u'_0(x)|\leq |b/\mu|+1$ for
$x\in [-b, -b+\delta]$, we then have $w(0, x)\geq u_0(x)$ in $[-b, b+\frac
1M]$ if $M\geq \max \{ \frac 1\delta, \frac
{|g_1|/\mu+1}{M_1}\}$. Making use of the comparison principle yields $u(t,
x)\leq w(t, x)$ in $\Omega$. Noticing that $u(t, g(t))=w(t,
g(t))=0$, we have
$$u_x(t, g(t))\leq w_x(t, g(t))=2M\overline M.$$
Recollecting the free boundary condition in (\ref{f1}) deduces
$$0<-g'(t)\leq 2\mu M\overline M\triangleq M_3,\quad 0<t\leq T,$$
where $M_3$ is independent of $T$.
Analogously, we can define
$$w(t, x)=\overline M[2M(h(t)-x)-M^2(h(t)-x)^2].$$
over the region
$$\Omega' =\{(t, x): \ 0<t\leq T,\  h(t)-\frac 1M<x<h(t)\}$$
get that $$0<h'(t)\leq M_4,\quad 0<t\leq T,$$
where $M_4$ is independent of $T$.
\end{pf}

\begin{thm} The solution of problem $(\ref{f1})$ exists and is
unique, and it can be extended to $[0, T_{\max})$ with
$T_{\max}\leq \infty$. Moreover, if $T_{\max}<\infty$, we have $\limsup_{t\to T_{\max}}
||u||_{L^\infty([g(t), h(t)]\times [0, t]}=\infty$.
\end{thm}
\begin{pf}
It follows from the uniqueness that there is a number $T_{\max}$ such
that $[0, T_{\max})$ is the maximal time interval in which the
solution exists. In order to prove the present theorem, it suffices
to show that, when $T_{\max}<\infty$, $\limsup_{t\to T_{\max}}
||u||_{L^\infty([0, t]\times[g(t), h(t)]}=\infty$. In what follows we
use the contradiction argument. Assume that $T_{\max}<\infty$ and
$||u||_{L^\infty([0, T_{\max})\times[g(t), h(t)])}<\infty$. Since $v\leq
M_2(M)$ in $[g(t), h(t)]\times [0, T_{\max})$ and $0<-g'(t)\leq M_3,0<h'(t)\leq M_4$ in
$[0, T_{\max})$ by Lemma 2.3, using a bootstrap argument and
Schauder's estimate yields a priori bound of
$||u(t,x)||_{C^{1+\alpha}[g(t),h(t)]}+||v(t,
x)||_{C^{1+\alpha}(-\infty,\infty)}$ for all $t\in [0, T_{\max})$. Let the
bound be $M_5$. It follows from the proof of Theorem 2.1 that there
exists a $\tau>0$ depending only on $\overline M$, $M_2, M_3,M_4$ and
$M_5$ such that the solution of problem (\ref{f1}) with the
initial time $T_{\max}-\tau/2$ can be extended uniquely to the time
$T_{\max}-\tau /2+\tau$ that contradicts the assumption. Thus the
proof is complete.
\end{pf}

\section{Global solution for the case $b_1 c_2>b_2 c_1$}
To obtain the global existence, we first derive a prior estimate for the solution of (\ref{f1}).

\begin{lem} If  $b_1 c_2>b_2 c_1$, then the solution of the free boundary problem $(\ref{f1})$ satisfies
$$0<u(t, x)\leq K_1\quad \textrm{for} \quad 0\leq t\leq T, \ g(t)< x< h(t),$$
$$0\leq v(x,t)\leq K_2\quad \textrm{for} \quad 0\leq t\leq T,\ -\infty< x<\infty,$$
where $K_i$ is independent of $T$ for $i=1,2$.
\end{lem}
\begin{pf}
Firstly we have that $u>0$ in $[g(t), h(t)]\times
[0, T]$ and $v\geq 0$ in $(-\infty, \infty)\times
[0, T]$ provided that solution exists.

Since the solution is classical in $[0, T]$,  there exists a  $\tilde K(T)$  such that
$u(t, x)\leq c_1\tilde K$ and $v(t, x)\leq \tilde K$.  Next we give the proof for $u(t, x)\leq K_1$ and $v(t, x)\leq K_2$,
where
$$ K_1 :=m \frac{a_1c_2+a_2c_1}{b_1c_2-b_2c_1} >\max_{[-b,b]} u_0(x), \quad
K_2 :=m \frac{a_1b_2+a_2b_1}{b_1c_2-b_2c_1} >||v_0||_{L^\infty(-\infty, \infty)} $$
for some $m>1.$

 Because  the interval $(-\infty, \infty)$ is unbounded,  maximum principle becomes invalid, next we prove that for any $l>b$,
  $$u(t, x)\leq K_1+\frac{(1+b_1)c_1}{b_1} \frac {\tilde K(x^2+2\tilde{d} t)}{l^2}, \quad
   v(t, x)\leq K_2+(1+b_1)\frac {\tilde K(x^2+2\tilde{d} t)}{l^2}$$
  for $0\leq t\leq T$, $ -l\leq x\leq l$, where $\tilde{d}=\max(d_1, d_2)$.
 Setting $$\ol u(t, x)= K_1+\frac{(1+b_1)c_1}{b_1} \frac {\tilde K(x^2+2\tilde{d} t)}{l^2},$$
 $$\ol v(t, x)=K_2+(1+b_1)\frac {\tilde K(x^2+2\tilde{d} t)}{l^2},$$
 then $(\ol u, \ol v)$ satisfies
 \begin{eqnarray*}
\left\{
\begin{array}{lll}
\ol u_{t}-d_1\ol u_{xx}\geq \ol u(a_1-b_{1}\ol u+c_{1}\ol v(t-\tau_1,x),\; & 0<t \leq T, \ -l<x<l,  \\
\ol v_{t}-d_2 \ol v_{xx}\geq \ol v(a_2+b_{2}\ol u(t-\tau_2,x)- c_{2}\ol v),\; & 0<t \leq T, \ -l<x<l, \\
\ol u \geq K_1+\frac{(1+b_1)c_1}{b_1} \tilde K>u,\ \ol v\geq K_2+(1+b_1)\tilde K>v,&0<t\leq T, \ x=\pm l, \\
\ol u(t, x)\geq K_1>u_{0}(x), & -\tau_2\leq t \leq 0,-l\leq x\leq l \\
\ol v(t, x)\geq K_2> v_{0}(x), &-\tau_1\leq t \leq 0,-l\leq x\leq l.
\end{array} \right.
\end{eqnarray*}
 It follows  that $ u\leq \ol u$ and $v\leq \ol v  $  by using the maximum principle on  $[0, T]\times [-l, l].$
Now for any fixed $(t_0, x_0)\in [0, T]\times (-\infty, \infty)$, let $l$ sufficiently large so that
$(t_0, x_0)\in [0, T]\times [-l, l]$, we deduce from the above proof that
 $$u(t_0, x_0)\leq \ol u(t_0, x_0)= K_1+\frac{(1+b_1)c_1}{b_1} \frac {\tilde K(x^2_0+2\tilde{d} t_0)}{l^2},$$
 $$v(t_0, x_0)\leq \ol v(t_0, x_0)=K_2+(1+b_1)\frac {\tilde K(x^2_0+2\tilde{d} t_0)}{l^2}.$$
Taking $l\to \infty$ gives the desired estimates.
\end{pf}

Combing Theorem 2.4 with Lemma 3.1 yields the following global existence:

\begin{thm} If parameters in double free boundaries problem $(\ref{f1})$ satisfy  $b_1 c_2>b_2 c_1$, then  $(\ref{f1})$ admits a unique global solution.
\end{thm}

Next we mainly give the long-time behavior of the free boundary problem (\ref{f1}). Here, we first give the slow solution.
\begin{thm} If $b_1 c_2>b_2 c_1$ and  $a_1>d_1(\frac{\pi}{2b})^2 $,
the free boundaries of the problem $(\ref{f1})$ satisfy
$$h_\infty=\infty\ \textrm{and} \ \ g_\infty=-\infty.$$
\end{thm}
\begin{pf} Combing Theorems 2.2 with Theorem  3.2, we know that the solution is global, $x=g(t)$ is
monotonic decreasing and $x=h(t)$ is monotonic increasing
Assume that $g_\infty
>-\infty$  by contradiction, then we have\\
  $\lim_{t\to +\infty}g'(t)=0.$

On the other hand, the condition $a_1>d_1(\frac{\pi}{2b})^2 $
implies that $a>\lambda_1$, where $\lambda_1$ denotes the first
eigenvalue of the problem
$$ -d_1\phi''=\lambda \phi \ \ \textrm{in}\ \ (-b, b),\quad \phi(\pm b)=0.$$
Therefore for all small $\delta>0$, the first eigenvalue
$\lambda^\delta _1$ of the problem
$$ -d_1\phi_{xx}-\delta \phi'=\lambda \phi \ \ \textrm{in}\ \ (-b, b),\quad \phi(\pm b)=0$$
satisfies $\lambda^\delta _1<a_1$. Fix such an $\delta >0$
and consider the problem
\begin{equation}
L_\delta \psi=a_1 \psi -b_1 \psi^2 \ \ \textrm{in}\ \ (-b,
b),\quad \psi(\pm b)=0,\label{m0}
\end{equation}
where $L_\delta \psi=-d_1 \psi ''-\delta \psi '$.  It is well
known (Proposition 3.3 in \cite{CC}) that the problem (\ref{m0})
admits a unique positive solution $\psi=\psi_\delta$. By the
moving plane method one easily sees that $\psi(x)$ is symmetric
about $x=0$ with  $\psi ' (x)< 0 $ for $x\in (0, b]$. Moreover using the
comparison principle, we have $\psi< \frac {a_1}{b_1}$ in $[-b, b]$.
We now set
$$ F(t, x)=\psi \left(\frac {b }{g(t)}x \right),$$
and directly compute
\begin{eqnarray*}
F_t -d_1 F_{xx}&=& \frac {-b x}{g^2(t)} g'(t) \psi '-d_1 \frac
{b^2}{g^2(t)} \psi ''\\
&=&\frac {b^2}{g^2(t)}[ -d_1 \psi ''+ \frac {x g'(t)}{-b} \psi '].
\end{eqnarray*}
Note that $g'(t)\to 0$ as $t\to +\infty$, we can choose $T_0>0$
such that $g'(t)>\delta \frac{b}{g_\infty}$ for $t\geq T_0$ and
hence for $t\geq T_0$ and $x\in [g(t), 0]$, we have
$\frac{xg'(t)}{-b}\geq -\delta.$ Therefore for such $t$ and
$x$,
\begin{eqnarray*}
F_t -d_1 F_{xx}&\leq & \frac {b^2}{g^2(t)}(-d_1 \psi ''-\delta
\psi ')\\
&=& \frac {b^2}{g^2(t)}(a_1 \psi -b_1 \psi^2).
\end{eqnarray*}
Because of  $0\leq \psi<\frac {a_1}{b_1}$, we have $a_1\psi -b_1
\psi^2\geq 0$ and hence from $\frac {-b}{g(t)}\leq 1$ we get
$$
F_t- d_1 F_{xx}\leq a_1\psi -b_1 \psi^2= a_1 F -b_1 F^2\quad
\textrm{for}\ \ t\geq T_0, \ x\in [g(t), 0].$$

 Now we choose $\varepsilon \in (0, 1)$ sufficient small so that $\varepsilon F(T_0,
x)\leq u(T_0, x)$. Then $\underline u(t, x):=\varepsilon F(t, x)$
 satisfies
\begin{eqnarray*}
\left\{
\begin{array}{lll}
\underline u_t-d_1 \underline u_{xx}\leq a_1 \underline u-b_1
\underline u^2,\; &t\geq T_0,\
x\in [g(t), 0], \\
\underline u(t, g(t))=0, \quad \underline u_x(t, 0)=0,  & t\geq T_0,\\
\underline u(T_0, x)\leq u(T_0, x), & 0\leq x\leq g(T_0).
\end{array} \right.
\end{eqnarray*}
So we can use the comparison principle to draw a conclusion
$$\underline u (t, x)\leq u(t, x)\ \ \textrm{for} \ t\geq T_0, \
x\in [g(t), 0].$$ It follows that
$$u_x(t, g(t))\geq \underline u_x(t, g(t)) =\delta \frac {b}{g(t)}
\psi '(b)\to \delta \frac {b}{g_\infty} \psi '(b)>0,$$ which means
that $g'(t)\leq -\mu \delta \frac {b}{g_\infty} \psi '(b)<0$. This
is a contradiction to the fact that $g'(t)\to 0$ as $t\to \infty$.
This contradiction implies that $g_\infty=-\infty.$
Likewise, we can set
$$F(t, x)=\psi \left(\frac {b }{h(t)}x \right), \ x\in [0, h(t)]$$
to prove that $ h_\infty=+\infty.$
\end{pf}

\section{Global and nonglobal  solutions}
In this section, we discuss the asymptotic behavior of the solution for the case $b_1 c_2<b_2 c_1$,
which is more intricate than that for the case $b_1 c_2>b_2 c_1$.
 First we present the blowup result.

\begin{thm} If $b_1 c_2<b_2 c_1$, then

$(i)$ the solution of  the free boundary problem $(\ref{f1})$  with any nontrival nonnegative
initial data blows up in case $a_i>d_i(\frac{\pi}{2b})^2 $ for $i=1,2$.

$(ii)$ the solution of  the free boundary problem $(\ref{f1})$ blows up for any $a_i$ in case the initial data
is sufficiently large .
\end{thm}
\begin{pf}

To prove this, it suffices to compare the free boundary problem with the corresponding problem in the fixed domain:
\begin{eqnarray}
\left\{
\begin{array}{lll}
u_{t}-d_1 u_{xx}=u(a_1-b_{1}u+c_{1}v(t-\tau_1,x)),\; &t>0,\, -b<x<b, \\
v_{t}-d_2v_{xx}=v(a_2+b_{2}u(t-\tau_2,x)-c_{2}v),\; & t>0, \, -b<x<b,\\
u(t, -b)=v(t, -b)=0,&t>0,\\
u(t, b)=v(t, b)=0,&t>0,\\
u(0, x)=u_{0}(x)\geq 0, &-\tau_2\leq t \leq 0,\,-b\leq x\leq b,\\
\ v(0, x)=v_{0}(x)\geq 0,&-\tau_1\leq t \leq 0,\,-b\leq x\leq b.
\end{array} \right.
\label{f13}
\end{eqnarray}
It follows from \cite{P} that the solution blows up if   $a_i>d_i(\frac{\pi}{2b})^2 $, $i=1,2$ or
if the initial data is sufficiently large.
We come to a conclusion  by making use of  maximum principle.
\end{pf}

\begin{rmk} The above theorem means that if the initial length $b$ is
large enough or if the initial data is sufficiently large, the solution will blow up. The constant $(\frac{\pi}{2b})^2 $
is the first eigenvalue of $-\Delta $ in $[-b, b]$ with homogeneous Dirichlet boundary condition.

\end{rmk}

The comparison principle used above is for the stationary boundary. In the following
we introduce a comparison principle for  double free boundaries $x=h(t)$ and $x=g(t)$.

\begin{lem} Suppose that $T\in (0,\infty)$, $\overline h,\overline g \in
C^1([0,T])$, $\overline u\in C(\overline D_{1,T}^*)\cap
C^{1,2}(D_{1,T}^*)$ and $\overline v\in C(\overline D_{2,T}^*)\cap
C^{1,2}(D_{2,T}^*)$ with $D_{1,T}^*=(0, T]\times (\overline g(t), \overline h(t))$,
$D_{2,T}^*=(0, T]\times (-\infty, +\infty)$, and
\begin{eqnarray*}
\left\{
\begin{array}{lll}
\overline u_{t}-d_1 \overline u_{xx}\geq \overline
u(a_1-b_{1}\overline u+c_{1}\overline v(t-\tau_1,x)),\;
& t>0, \, \overline g(t)<x<\overline h(t),  \\
\overline v_{t}-d_2 \overline v_{xx}\geq \overline
v(a_2+b_{2}\overline u(t-\tau_2,x)-c_{2}\overline v),\;
& t>0, \, -\infty<x<\infty, \\
\overline u(t, x)=0,\; & t>0, \  -\infty<x<g(t),  \\
\overline u(t, x)=0,\; & t>0, \  h(t)<x<\infty,  \\
\overline u=0,\quad \overline h'(t)\geq -\mu \frac{\partial
\overline u}{\partial
x},\quad & t>0, \  x=\overline h(t),  \\
\overline u=0,\quad \overline g'(t)\leq -\mu \frac{\partial
\overline u}{\partial
x},\quad & t>0, \  x=\overline g(t).
\end{array} \right.
\end{eqnarray*}
If $-b\geq \overline g(0), b\leq \overline h(0)$, $u_0(x)\leq \overline u(t,x)$
in $[-b,b]\times[-\tau_2,0]$ and $v_0(x)\leq \overline v(t,x)$
in $(-\infty,+\infty)\times[-\tau_1,0]$, then the solution $(u, v,g, h)$ of the free boundary problem $\eqref{f1}$
satisfies
$$g(t)\geq\overline g(t),\ h(t)\leq\overline h(t)\ \textrm{in}\ (0, T],$$
$$\ u(t,x)\leq\overline u(t,x)\  \textrm{in}\ [0, T]\times (g(t), h(t))$$
$$\ v(t,x)\leq \overline v(t,x)\ \textrm{in}\ [0, T]\times (-\infty,+\infty).$$
\end{lem}
\begin{pf} We first suppose that $g(0)>\overline g(0),h(0)<\overline h(0)$. In this case
we assert that $g(t)>\overline g(t)$ and $h(t)<\overline h(t)$ for $0<t\leq T$ by using  contradiction.
If it is not true, then there exists $t^*_1\in (0, T]$
 such that $h(t)<\overline h(t)$ for $t\in [0, t^*_1)$ and $h(t^*_1)=\overline h(t^*_1)$.
It follows that
$$h'(t^*_1)\geq \overline h'(t^*_1).$$
Because of  the system of (\ref{f1}) is nondecreasing, and
applying the strong maximum principle for the parabolic systems give that
$u(t,x)<\overline u(t,x)$ in $(0, t^*_1]\times (g(t), h(t))$ and $$\frac {\partial}{\partial x}(u-\overline u)|_{(t^*_1, h(t^*_1))}>0$$
by $u(t^*_1, h(t^*_1))=0=\overline u(t^*_1, \overline h(t^*_1))$,
then $$h'(t^*_1)=-\mu \frac{\partial u}{\partial x}(t^*_1,
h(t^*_1))<\overline h'(t^*_1).$$ This leads to a contradiction,
which proves our assert that $h(t)<\overline h(t)$ for $0<t\leq T$ if
$h(0)=b<\overline h(0)$. Analogously, we can prove that $g(t)>\overline g(t) $ for $ 0<t\leq T.$
 Now we may draw a conclusion that
  $u(t,x)\leq \overline u(t,x)$ in $[0, T]\times (g(t), h(t))$ and
$v(t,x)\leq \overline v(t,x)$ in $[0, T]\times (-\infty,+\infty)$ by approximation.
 \end{pf}

\begin{rmk} The  $(\overline u, \overline v, \overline h,\overline g)$ in Lemma 4.2 is
usually called an upper solution of the problem \eqref{f1}. We can
define a lower solution by reversing all the  inequalities in the
obvious places. Moreover, one can easily prove an analogue of Lemma
4.2 for lower solutions.
\end{rmk}

 In the following theorem ,we show existence of a global fast solution .

\begin{thm} If $b_1 c_2<b_2 c_1$, then the free boundary problem $(\ref{f1})$
admits a global fast solution provided that the initial data $u_0$ and $b$ are suitably small.
Moreover, there exist constant $C,\beta >0$ depending on $b, u_0 $ and $k$
such that
$$||u||_\infty\leq C e^{-\beta t},\quad t\geq 0 $$
for some $k>1$.
\end{thm}
\begin{pf}
  \cite{RT}was the main source of inspiration for its proof, we have only to the structure proper global supersolution.
  Define
$$\sigma (t)=2b(k-e^{-\gamma t}), \lambda (t)=-\sigma (t), \  t\geq 0, \quad W(y)=\cos (\frac{\pi}{2}y), \ -1\leq y\leq 1,$$
and $$\overline u(t, x)=\delta e^{-\beta t}W(x/\sigma (t)), \ t\geq 0,\ \lambda(t)\leq x\leq \sigma(t).$$
$$\overline v(t, x)=k\frac {a_2}{c_2},\ t\geq 0,\ -\infty\leq x\leq \infty,$$
where $\gamma, \beta$ and $\delta>0$ to be determined later. \\
 Straightforward calculations yields
\begin{eqnarray*}
& &\overline u_t-d_1\overline u_{xx}-\overline u(a_1-b_1\overline u+c_1\overline v(t-\tau_1,x))\\
& &=\delta e^{-\beta t}[-\beta W-x\sigma '\sigma^{-2}W'-d_1\sigma^{-2}W''-W(a_1-
b_1\delta e^{-\beta t}W+c_1\overline v)]\\
& &\geq \delta e^{-\beta t}W[-\beta +(\frac \pi 2)^2\frac {d_1}{4k^2b^2}-a_1-k\frac {c_1a_2}{c_2}]
\end{eqnarray*}
for all $t>0$  and  $\lambda (t)<x<\sigma (t)$ and
\begin{eqnarray*}
& &\overline v_{t}-d_2\overline v_{xx}-\overline v(a_2+b_{2}\overline u(t-\tau_2,x)-c_{2}\overline v)\\
& &=k\frac {a_2}{c_2}(-a_2-b_2\delta e^{-\beta (t-\tau_2)} W+ka_2)\geq k\frac {a_2}{c_2}((k-1)a_2-b_2\delta e^{\beta\tau_2})
\end{eqnarray*}
for all $t>0$ and $-\infty<x<\infty$. On the other hand, we can easily deduce
$\sigma'(t)=2\gamma be^{-\gamma t}>0$ , $-\overline u_x(t, \sigma
(t))=\frac {\pi}{2}\delta \sigma^{-1}(t)e^{-\beta t}$ and $-\overline u_x(t, \lambda
(t))=\frac {\pi}{2}\delta \lambda^{-1}(t)e^{-\beta t}$. Now we set $b_0$
such that
$$\frac {d_1}{8k^2b_0^2}(\frac \pi 2)^2=a_1+\frac {ka_2c_1}{c_2},$$
if $0<b\leq b_0$, setting
$$\delta =\min \{\frac {(k-1)a_2}{b_2e^{\beta\tau_2}}, \ \frac {(k-1)d_1\pi}{2k^2\mu}(\frac { b} {2b_0})^2\}, \
\beta=\gamma=(\frac \pi 2)^2\frac{d_1}{16k^2b^2_0},$$
 It follows that
\begin{eqnarray*}
\left\{
\begin{array}{lll}
\overline u_{t}-d_1 \overline u_{xx}\geq \overline u(a_1-b_{1}\overline u+c_{1}\overline v(t-\tau_1,x)),\; &t>0,\ \lambda<x<\sigma(t),  \\
\overline v_{t}-d_2 \overline v_{xx}\geq \overline v(a_2+b_{2}\overline(t-\tau_2,x)-c_{2}\overline v),\; &t>0, \  -\infty<x<\infty, \\
\overline u=0,\quad \sigma'(t)> -\mu \frac{\partial \overline u}{\partial
x},\quad &t>0, \  x=\sigma(t),\\
\overline u=0,\quad \lambda'(t)< -\mu \frac{\partial \overline u}{\partial
x},\quad &t>0, \  x=\lambda(t),\\
\sigma (0)=2b>b, \lambda (0)=-2b<-b.&
\end{array} \right.
\end{eqnarray*}
By making use of the maximum principle, we can get that $h(t)<\sigma(t),$  $g(t)>\lambda(t),$
and $u(t, x)<\overline u(t, x)$, $v(t, x)<\overline v(t, x)$ for $g(t)\leq x\leq h(t)$ provided
 $(u,v)$ exists.  particularly, it follows from Lemma
4.2 that $(u, v)$ exists globally and $
g_\infty>-\infty$,   $
h_\infty<\infty$.
\end{pf}

From the above proof, we have the following global existence result
\begin{thm} If $b_1 c_2<b_2 c_1$ and  $a_1\leq 0, a_2\leq 0$, then  the free boundary problem $(\ref{f1})$
admits a global fast solution provided $u_0$ is suitably small.
\end{thm}

\begin{rmk} If $b_1 c_2>b_2 c_1$, Theorem 3.3
shows that the solution is slow for any initial data. If $b_1
c_2<b_2 c_1$, Theorem 4.1 shows that the solution blows up for large initial data, and sufficient conditions for the global fast solution are given
in Theorems 4.3 and 4.4, which implies that the global fast solution is possible if the initial data is suitably small.
\end{rmk}

\end{document}